\documentclass[12pt]{article}
\usepackage{ucs}
\usepackage[utf8x]{inputenc}
\usepackage{amsfonts}
\usepackage{amsmath}
\usepackage[english]{babel}
\usepackage{indentfirst}
\usepackage{geometry}
\usepackage{graphicx}
\usepackage{multicol}

\title{Diophantine approximations with Fibonacci numbers}
\author{Victoria Zhuravleva}
\date{}

\begin{document}

\maketitle

\begin{abstract}

Let $F_{n}$ be the $n$-th Fibonacci number. Put $\varphi=\frac{1+\sqrt5}{2}$. We prove that the following inequalities hold for any real $\alpha$:

1) $\inf_{n \in \mathbb N } ||F_n\alpha||\le\frac{\varphi-1}{\varphi+2}$,

2) $\liminf_{n\to \infty}||F_n\alpha||\le \frac{1}{5}$,

3) $\liminf_{n \to \infty}||\varphi^n \alpha||\le \frac{1}{5}$.

These results are the best possible.

\end{abstract}

\section{Introduction}

In this paper $||\alpha||$ denotes the distance from a real $\alpha$ to the nearest integer. Let 
$$F_1=F_2=1, F_3=2, F_4=3, F_5=5, F_6=8, F_7=13,\ldots$$ be Fibonacci numbers. For convenience we put $F_{-1}=1$, $F_0=0$.

As $F_{n}$ are distinct integers we deduce from H. Weyl's Theorem (see \cite{ref4}, Ch.1 \S4) that the fractional parts $\{F_n\alpha\}$ are uniformly distributed for almost all real $\alpha$. From the other hand Fibonacci numbers form a lacunary sequence. So the set 

$$
\mathcal N=\{\alpha: \exists\gamma(\alpha)>0 \text{ such that } \inf_{n \in \mathbb N} ||F_n \alpha||\ge\gamma(\alpha)  \}
$$
is an $\alpha$-winning set for every $\alpha \in (0, 1/2]$ in the sence of W.M. Schmidt's $(\alpha, \beta)$-games and hence the Hausdorff dimension of the set $\mathcal N$ is equal to one. For the definition and simplest properties of winning sets see Ch.3 from \cite{ref5} and \cite{ref6}. Some quantitative version of W.M. Schmidt's results one can find in \cite{ref1}.

Put $\varphi=\frac{1+\sqrt5}{2}$. As 
\begin{equation}
F_{n+1}=F_{n}+F_{n-1}, \quad \varphi^{n+1}=\varphi^n+\varphi^{n-1}
\end{equation}
one can easily see that for any real $\alpha$ and for any positive integer $n$ the following inequlities are valid:

\begin{equation}
\min_{j=n-1, n, n+1} ||\varphi^j \alpha||\le\frac13, \quad \min_{j=n-1, n, n+1} ||F_{j}\alpha||\le\frac13
\end{equation}

\section{A result by A. Dubickas}

In \cite{ref3} A. Dubickas proved a result related to Diophantine approximations with powers of algebraic numbers. Here we would like to give the formulation of this result. 

The length $L(P)$ of a polynomial $P(x)=p_0+p_1x+\ldots+p_kx^k \in \mathbb R [x]$ is defined as the sum of absolute values of all coefficients of $P(x)$:
$$
L(P)=|p_0|+|p_1|+\ldots+|p_k|.
$$
The \textit{reduced length} $l(P)$ of a polynomial $P(x)$ is defined as
$$
l(P)=\inf_{Q}L(PQ),
$$
where the infinum is taken over all polynomials $Q(x)=q_0+q_1x+\ldots+q_rx^r \in \mathbb R [x]$ such that $q_0=1$ or $q_r=1$. The reduced length $l(\alpha)$ of an algebraic number $\alpha$ is defined as the reduced length of the irreducible polynomial $P_{\alpha}(x)\in \mathbb Z [x]$ such that $P_{\alpha}(\alpha)=0$.

\textsc{Theorem}[A. Dubickas, \cite{ref3}].
Suppose $\tau>1$ is an arbitrary algebraic number. Suppose $\alpha$ be a positive real number that lies outside the field $\mathbb Q (\tau)$ if $\tau$ is a Pisot or a Salem number. Then it is not possible that all fractional parts of the form $\{\tau^j \alpha\}$, $j \in \mathbb N$ belong to a certain open interval of the length $1/l(\tau)$.

The history of the question as well as the definitions of Pisot and Salem numbers one can find in \cite{ref3}. Here we should note that $\varphi=\frac{1+\sqrt5}{2}$ is a Pisot number and $l(\varphi)=1+\varphi$ (see \cite{ref3}). Particularly A. Dubickas's theorem shows that under the condition $\alpha \in \mathbb R \backslash \mathbb Q (\varphi)$ for any $n_0$ there exists an integer $n\ge n_0$ such that
$$
\{\varphi^n\alpha\} \notin \Bigr (\frac12-\frac{1}{2(1+\varphi)}, \frac12+\frac{1}{2(1+\varphi)} \Bigr)
$$ 
So A.Dubickas's theorem leads to the following asymtotic inequality:
\begin{equation}
\liminf_{n \to \infty} ||\varphi^n \alpha||\le\frac{1}{2\varphi}=\frac{\sqrt 5 - 1}{4}.
\end{equation}
As for Fibonacci numbers we have the formula
$$
F_n=\frac{1}{\sqrt5}\Bigr (\varphi^n-\Bigr (-\frac{1}{\varphi}\Bigr)^n\Bigr)
$$
we immediately deduce

\begin{equation}
\liminf_{n \to \infty} ||F_n \alpha||\le\frac{1}{2\varphi}=\frac{\sqrt 5 - 1}{4}
\end{equation}
in the case $\alpha \in \mathbb R \backslash \mathbb Q (\varphi)$.

\section{Statement of results}

In this section we summarize all results obtained in this paper. 

Let $K$, $N$ be positive integers. Put $$d_N^K=\max_{\alpha  \in \mathbb R} \min_{k=K, \ldots ,K+N-1} ||F_k \alpha||.$$

\textsc{Theorem 1}.

1) \textit{The following equalities are valid: $d_1^1=d_2^1=\frac{1}{2}$, $d_3^1=\frac{1}{3}$, $d_4^1=d_5^1=\frac{1}{4}$.}

2) \textit{Let $N \ge 6$, put $n=[\frac{N-2}{4}]$. Then $d^1_{N}=\frac{F_{2n+1}}{F_{2n+2}+F_{2n+4}}$.}
\\

\textsc{Corollary 1}. \textit{The following equality is valid: }$$\lim_{N \to \infty} d^1_N=\frac{\varphi-1}{\varphi+2}.$$

\textsc{Corrollary 2}. \textit{The following inequality holds for any real $\alpha$:}
$$\inf_{n \in \mathbb N}||F_n \alpha||\le \frac{\varphi-1}{\varphi+2}.$$

\textsc{Theorem 2.} \textit{For $N\ge4$, $\alpha_1=\frac{1}{\varphi+2}$ one has}
$$\min_{n \le N}||F_{n}\alpha_1||=\frac{\varphi-1}{\varphi+2}.$$

\textsc{Theorem 3.} \textit{Let $\alpha_1=\frac{1}{\varphi+2}$. Then $\forall \varepsilon >0$ $\forall N$ $\exists K=K(N)$: $$\min_{k=K \ldots K+N}||F_k \alpha_1||>\frac{1}{5}-\varepsilon.$$}

\textsc{Theorem 4.} \textit{Let $G_1$ and $G_2$ be arbitrary real numbers, and $G_n=G_{n-1}+G_{n-2}$ for $n\ge3$. Then:}

1) $\max_{G_1, G_2} \min_{n=1,2}||G_n||= \frac12$,

2) $\max_{G_1, G_2} \min_{n=1,2,3}||G_n||= \frac13$,

3) $\max_{G_1, G_2} \min_{n=1,2,3,4}||G_n||= \frac14$,

4) $\max_{G_1, G_2} \min_{n=1,2,3,4,5}||G_n||= \frac14$,

5) $\max_{G_1, G_2} \min_{n=1, \ldots ,k}||G_n||= \frac15$ for $k\ge 6$.
\\

\textsc{Corollary 3.} \textit{
Suppose $N\ge 6$, then }$$\lim_{K\to \infty} d_N^K=\frac{1}{5}.$$

\textsc{Corollary 4.} \textit{Suppose $N\ge 6$, then}
$$ \lim_{K\to \infty} \max_{\alpha \in \mathbb R} \min_{k=K,\ldots,K+N-1}||\varphi^k \alpha||=\frac{1}{5}. $$

\textsc{Corollary 5.} \textit{The following equalities are valid for any real $\alpha$:} $$\liminf_{n \to \infty}||F_n \alpha||\le\frac{1}{5},\quad \liminf_{n\to\infty}||\varphi^n \alpha||\le\frac{1}{5}.$$

One can see that we improve inequlities (3) and (4). So we improve A. Dubickas's Theorem, but only for $\tau=\varphi$.

\section{The proof of Theorem 1 for $N<6$}

The function $||F_kx||$ is periodic with period less than or equal to 1. Also it is symmetric with the respect to the line $x=1/2$. Thus without loss of generality we consider this function on the segment $[0, 1/2]$. 

Since $||F_kx||$ is a piecewise linear function its graph consists of line segments. Let $t$ be an arbitrary integer. Then

$$
||F_kx||=\begin{cases}
t-F_kx &\text{if $x \in [\frac{t}{F_k}-\frac{1}{2F_k}, \frac{t}{F_k}]$,}\\
F_kx-t &\text{if $x \in [\frac{t}{F_k}, \frac{t}{F_k}+\frac{1}{2F_k}]$.}
\end{cases}
$$

Put $F_N(x)=\min_{k=1,\ldots,N}||F_kx||.$
\\

\textsc{Lemma 1.} \textit {One can easily see that for $x \in [0, \frac{1}{2}]$}
\begin{multicols}{2}
$$
F_1(x)=F_2(x)=x,
$$
$$
F_3(x)=\begin{cases}
x &\text{if $x \in [0,\frac{1}{3}]$,}\\
1-2x &\text{if $x \in [\frac{1}{3}, \frac{1}{2}]$,}
\end{cases}
$$
$$
F_4(x)=\begin{cases}
x &\text{if $x \in [0, \frac{1}{4}]$,}\\
1-3x &\text{if $x \in [\frac{1}{4}, \frac{1}{3}]$,}\\
3x-1 &\text{if  $x \in [\frac{1}{3}, \frac{2}{5}]$,}\\
1-2x &\text{if $x \in [\frac{2}{5}, \frac{1}{2}]$,}
\end{cases}
$$
$$
F_5(x)=\begin{cases}
x &\text{if $x \in [0, \frac{1}{6}]$,}\\
1-5x &\text{if $x \in [\frac{1}{6}, \frac{1}{5}]$,}\\
5x-1 &\text{if $x \in [\frac{1}{5}, \frac{1}{4}]$,}\\
1-3x &\text{if $x \in [\frac{1}{4}, \frac{1}{3}]$,}\\
3x-1 &\text{if $x \in [\frac{1}{3}, \frac{3}{8}]$,}\\
2-5x &\text{if $x \in [\frac{3}{8}, \frac{2}{5}]$,}\\
5x-2 &\text{if $x \in [\frac{2}{5}, \frac{3}{7}]$,}\\
1-2x &\text{if $x \in [\frac{3}{7}, \frac{1}{2}]$,}
\end{cases}
$$
\\
$$
F_6(x)=\begin{cases}
x &\text{if $x \in [0, \frac{1}{9}]$,}\\
1-8x &\text{if $x \in [\frac{1}{9},\frac{1}{8}]$,}\\
8x-1 &\text{if  $x \in [\frac{1}{8}, \frac{1}{7}]$,}\\
x &\text{if $x \in [\frac{1}{7}, \frac{1}{6}]$,}\\
1-5x &\text{if $x \in [\frac{1}{6}, \frac{1}{5}]$,}\\
5x-1 &\text{if $x \in [\frac{1}{5}, \frac{3}{13}]$,}\\
2-8x &\text{if $x \in [\frac{3}{13}, \frac{1}{4}]$,}\\
8x-2 &\text{if $x \in [\frac{1}{4}, \frac{3}{11}]$,}\\
1-3x &\text{if $x \in [\frac{3}{11}, \frac{1}{3}]$,}\\
3x-1 &\text{if $x \in [\frac{1}{3}, \frac{4}{11}]$,}\\
3-8x &\text{if $x \in [\frac{4}{11}, \frac{3}{8}]$,}\\
8x-3 &\text{if $x \in [\frac{3}{8}, \frac{5}{13}]$,}\\
2-5x &\text{if $x \in [\frac{5}{13}, \frac{2}{5}]$,}\\
5x-2 &\text{if $x \in [\frac{2}{5}, \frac{3}{7}]$,}\\
1-2x &\text{if $x \in [\frac{3}{7}, \frac{1}{2}]$,}
\end{cases}
$$
\end{multicols}

\newpage
\begin{multicols}{2}
$$
F_7(x)=\begin{cases}
x, &\text{for $x \in [0, \frac{1}{14}]$}\\
1-13x, &\text{for $x \in [\frac{1}{14}, \frac{1}{13}]$}\\
13x-1, &\text{for $x \in [\frac{1}{13}, \frac{1}{12}]$}\\
x, &\text{for $x \in [\frac{1}{12}, \frac{1}{9}]$}\\
1-8x, &\text{for $x \in [\frac{1}{9},\frac{1}{8}]$}\\
8x-1, &\text{for  $x \in [\frac{1}{8}, \frac{1}{7}]$}\\
2-13x, &\text{for $x \in [\frac{1}{7}, \frac{2}{13}]$}\\
13x-2, &\text{for $x \in [\frac{2}{13}, \frac{1}{6}]$}\\
1-5x, &\text{for $x \in [\frac{1}{6}, \frac{1}{5}]$}\\
5x-1, &\text{for $x \in [\frac{1}{5}, \frac{2}{9}]$}\\
3-13x, &\text{for $x \in [\frac{2}{9}, \frac{3}{13}]$}\\
13x-3, &\text{for $x \in [\frac{3}{13}, \frac{4}{17}]$}\\
2-8x, &\text{for $x \in [\frac{4}{17}, \frac{1}{4}]$}\\
8x-2, &\text{for $x \in [\frac{1}{4}, \frac{3}{11}]$}\\
1-3x, &\text{for $x \in [\frac{3}{11}. \frac{3}{10}]$}
\end{cases}
$$

$$
F_7(x)=\begin{cases}
4-13x &\text{if $x \in [\frac{3}{10}, \frac{4}{13}]$,}\\
13x-4 &\text{if $x \in [\frac{4}{13}, \frac{5}{16}]$,}\\
1-3x &\text{if $x \in [\frac{5}{16}, \frac{1}{3}]$,}\\
3x-1 &\text{if $x \in [\frac{1}{3}, \frac{4}{11}]$,}\\
3-8x &\text{if $x \in [\frac{4}{11}, \frac{3}{8}]$,}\\
8x-3 &\text{if $x \in [\frac{3}{8}, \frac{8}{21}]$,}\\
5-13x &\text{if $x \in [\frac{8}{21}, \frac{5}{13}]$,}\\
13x-5 &\text{if $x \in [\frac{5}{13}, \frac{7}{18}]$,}\\
2-5x &\text{if $x \in [\frac{7}{18}, \frac{2}{5}]$,}\\
5x-2 &\text{if $x \in [\frac{2}{5}, \frac{3}{7}]$,}\\
1-2x &\text{if $x \in [\frac{3}{7}, \frac{5}{11}]$,}\\
6-13x &\text{if $x \in [\frac{5}{11}, \frac{6}{13}]$,}\\
13x-6 &\text{if $x \in [\frac{6}{13}, \frac{7}{15}]$,}\\
1-2x &\text{if $x \in [\frac{7}{15}, \frac{1}{2}]$.}
\end{cases}
$$
\end{multicols}

Put $d_N=\max_{x \in \mathbb R} F_N(x)$. Let $x_N$ be the point where the function $F_N(x)$ attains its maximal value. 

From Lemma 1 we find $x_N$ and $d_N$ for $N=1,\ldots,7$ (see Table 1). So we obtain the first statement of Theorem 1.
\\

\textsc{Lemma 2.} \textit {The graph of the function $F_{7}(x)$ (see Fig.1) has only one vertex which lies above the line $y=\frac{\varphi-1}{\varphi+2}$. This vertex has coordinates $(\frac{3}{11}, \frac{2}{11})$.}
\\

\newpage

\begin{figure}[h]
\includegraphics[width=\textwidth]{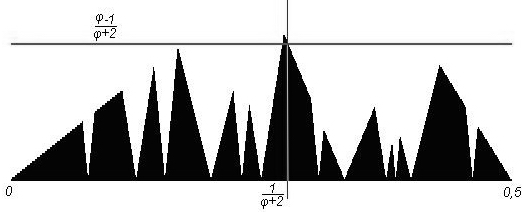}
\label{labelname}
\caption{$F_7(x)$}
\end{figure}

We continue to calculate the values of $x_N$ and $d_N$ (see Table 1).

\begin{center}
\begin{tabular}{|c|c|c|}
\hline
$N$ & $x_N$ & $d_N$\\
\hline
\ttfamily $1,2$ & $1/2$ & $1/2$\\
\ttfamily $3$ & $1/3$ & $1/3$\\
\ttfamily $4,5$ & $1/4$ & $1/4$\\
\ttfamily $6 \ldots 9$ & $3/11$ & $2/11$\\
\ttfamily $10 \ldots 13$ & $8/29$ & $5/29$\\
\ttfamily $14 \ldots 17$ & $21/76$ & $13/76$\\
\ttfamily $18$ & $55/199$ & $34/199$\\
\hline
\end{tabular}

\quad

Table 1.
\end{center}

From these results we note that
\[
x_6=\frac{3}{11}=\frac{F_4}{F_4+F_6}, \text{ }  d_6=\frac{2}{11}=\frac{F_3}{F_4+F_6},
\]
\[
x_{10}=\frac{8}{29}=\frac{F_6}{F_6+F_8}, \text{ }  d_{10}=\frac{5}{29}=\frac{F_5}{F_6+F_8},
\]
\[
x_{14}=\frac{21}{76}=\frac{F_8}{F_8+F_{10}}, \text{ }  d_{14}=\frac{13}{76}=\frac{F_7}{F_8+F_{10}},
\]
\[
x_{18}=\frac{55}{199}=\frac{F_{10}}{F_{10}+F_{12}}, \text{ }  d_{18}=\frac{34}{199}=\frac{F_{9}}{F_{10}+F_{12}}.
\]

One can see that points $(3/11, 2/11)$, $(8/29, 5/29)$, $(21/76, 13/76)$, $(55/199, 34/199)$ are the elements of the sequence of the points $(a_n, b_n)$ where $a_n=\frac{F_{2n+2}}{F_{2n+2}+F_{2n+4}}$, $b_n=\frac{F_{2n+1}}{F_{2n+2}+F_{2n+4}}$.

We note that $\lim_{n \to \infty} a_n=\frac{1}{\varphi+2}$ and $\lim_{n \to \infty} b_n=\frac{\varphi-1}{\varphi+2}$.

\section{The nearest integer to $\frac{F_n}{\varphi+2}$}

Let $T_n=\sum_{k=1}^{[n/2]}(-1)^{k+1}F_{n-2k}=F_{n-2}-F_{n-4}+F_{n-6}-F_{n-8}+\ldots$
\\

\textsc{Lemma 3.} \textit{$T_n$ is the nearest integer to $\frac{F_n}{\varphi+2}$.}

Proof.

Consider sums

$$S_{4t}=F_0+F_4+F_8+\ldots+F_{4t},$$ 
$$S_{4t+1}=F_1+F_5+F_9+\ldots+F_{4t+1},$$ 
$$S_{4t+2}=F_2+F_6+F_{10}+\ldots+F_{4t+2},$$ 
$$S_{4t+3}=F_3+F_7+F_{11}+\ldots+F_{4t+3}.$$

From (1) one can see that these sums sutisfy the following system of linear equations:

\[
\begin{cases}
S_{4t}+S_{4t+1}+S_{4t+2}+S_{4t+3}=F_{4t+5}-1,\\
S_{4t}+S_{4t+1}=S_{4t+2},\\
S_{4t+1}+S_{4t+2}=S_{4t+3},\\
S_{4t+3}+S_{4t}=S_{4t+1}+F_{4t+3}-F_1.
\end{cases}
\]

Therefore

$$S_{4t}=\frac{(4F_{4t+3}-F_{4t+5}-3)}{5},$$
$$S_{4t+1}=\frac{(2F_{4t+5}-3F_{4t+3}+1)}{5},$$
$$S_{4t+2}=\frac{(F_{4t+5}+F_{4t+3}-2)}{5},$$
$$S_{4t+3}=\frac{(3F_{4t+5}-2F_{4t+3}-1)}{5}.$$

Then the explicit formulas for $T_n$ are

$$T_{4t}=S_{4(t-1)+2}-S_{4(t-1)}=\frac{2F_{4t+1}-3F_{4t-1}+1}{5},$$

$$T_{4t+1}=S_{4(t-1)+3}-S_{4(t-1)+1}=\frac{F_{4t+1}+F_{4t-1}-2}{5},$$

$$T_{4t+2}=S_{4t}-S_{4(t-1)+2}=\frac{4F_{4t+3}-F_{4t+5}-F_{4t+1}-F_{4t-1}-1}{5},$$

$$T_{4t+3}=S_{4t+1}-S_{4(t-1)+3}=\frac{2F_{4t+5}-3F_{4t+3}-3F_{4t+1}+2F_{4t-1}+2}{5}.$$

For convenience we rewrite these expressions in the following form:

\begin{equation}
\label{eq2} T_{4t}=\frac{2F_{4t}-F_{4t-1}+1}{5},
\end{equation}

\begin{equation}
\label{eq3} T_{4t+1}=\frac{F_{4t}+2F_{4t-1}-2}{5},
\end{equation}

\begin{equation}
\label{eq4} T_{4t+2}=\frac{3F_{4t}+F_{4t-1}-1}{5},
\end{equation}

\begin{equation}
\label{eq5} T_{4t+3}=\frac{4F_{4t}+3F_{4t-1}+2}{5}.
\end{equation}

No we calculate the difference $T_n-\frac{F_n}{\varphi+2}$ using Binet's formula.

1) If $n=4t$ then
$$T_{4t}-\frac{F_{4t}}{\varphi+2} =\frac{1}{5}+\frac{2\varphi^{4t}-2(1-\varphi)^{4t}-\varphi^{4t-1}+(1-\varphi)^{4t-1}}{5\sqrt5}-\frac{\varphi^{4t}-(1-\varphi)^{4t}}{\sqrt5(\varphi+2)}=$$

$$=\frac{1}{5}+\frac{2\varphi^{4t+1}+4\varphi^{4t}-\varphi^{4t}-2\varphi^{4t-1}-5\varphi^{4t}}{5\sqrt5(\varphi+2)}+$$

$$+\frac{\varphi(1-\varphi)^{4t-1}-2\varphi(1-\varphi)^{4t}+2(1-\varphi)^{4t-1}-4(1-\varphi)^{4t}+5(1-\varphi)^{4t}}{5\sqrt5(\varphi+2)}=$$

$$=\frac{1}{5}+\frac{\varphi^{4t-1}(2\varphi^2-2\varphi-2)}{5\sqrt5(\varphi+2)}+\frac{(1-\varphi)^{4t-2}((1-\varphi)^2+4(1-\varphi)-1)}{5\sqrt5(\varphi+2)}=$$

$$=\frac{1}{5}+\frac{(1-\varphi)(1-\varphi)^{4t-2}}{\sqrt5(\varphi+2)}.$$

Similar formulas are obtained in the three remaining cases.

2) If $n=4t+1$ then
$$T_{4t+1}-\frac{F_{4t+1}}{\varphi+2}=-\frac{2}{5}+\frac{(2-\varphi)(1-\varphi)^{4t-2}}{\sqrt5(\varphi+2)}.$$

3) If $n=4t+2$ then
$$T_{4t+2}-\frac{F_{4t+2}}{\varphi+2}=-\frac{1}{5}+\frac{(3-2\varphi)(1-\varphi)^{4t-2}}{\sqrt5(\varphi+2)}.$$

4) If $n=4t+3$ then
$$T_{4t+3}-\frac{F_{4t+3}}{\varphi+2}=\frac{2}{5}+\frac{(5-3\varphi)(1-\varphi)^{4t-2}}{\sqrt5(\varphi+2)}.$$

For convenience we define the function $r(n)$ so that the following formulas are valid:

\begin{equation}
\label{eq6} T_{4t}-\frac{F_{4t}}{\varphi+2}=\frac{1}{5}+r(4t),
\end{equation}

\begin{equation}
\label{eq7} T_{4t+1}-\frac{F_{4t+1}}{\varphi+2}=-\frac{2}{5}+r(4t+1),
\end{equation}

\begin{equation}
\label{eq8} T_{4t+2}-\frac{F_{4t+2}}{\varphi+2}=-\frac{1}{5}+r(4t+2),
\end{equation}

\begin{equation}
\label{eq9} T_{4t+3}-\frac{F_{4t+3}}{\varphi+2}=\frac{2}{5}+r(4t+3).
\end{equation}

We see that $|r(4t)|>|r(4t+1)|>|r(4t+2)|>|r(4t+3)|>|r(4t+4)|$. So the function $|r(n)|$ is decreasing. 

To complete the proof of Lemma 3 we need to bound $r(4t)$, $r(4t+1)$, $r(4t+2)$ and $r(4t+3)$ for $t=1$. 

We use the bounds $1,618<\varphi<1,619$ and $2,236<\sqrt5<2,237$. Then 

\begin{equation}
\label{eq10} -0,030<\frac{(1-\varphi)(1-\varphi)^{2}}{\sqrt5(\varphi+2)}=r(4)<-0,029,
\end{equation}

\begin{equation}
\label{eq11} 0,018<\frac{(2-\varphi)(1-\varphi)^2}{\sqrt5(\varphi+2)}=r(5)<0,019,
\end{equation}

\begin{equation}
\label{eq12} -0,012<\frac{(3-2\varphi)(1-\varphi)^2}{\sqrt5(\varphi+2)}=r(6)<-0,011,
\end{equation}

\begin{equation}
\label{eq13} 0,006<\frac{(5-3\varphi)(1-\varphi)^2}{\sqrt5(\varphi+2)}=r(7)<0,007.
\end{equation}

So we proved lemma for $t \ge 1$. For $t=0$ Lemma 3 can be verified  directly: $0$ is the nearest integer to $\frac{F_{1}}{\varphi+2}=\frac{F_{2}}{\varphi+2}=\frac{1}{\varphi+2}$, $1$ - to $\frac{F_{3}}{\varphi+2}=\frac{2}{\varphi+2}$.

Lemma 3 is proved.

\section{Proof of Theorem 1}

We fix an arbitrary integer $t\ge 1$.

To prove the second part of Theorem 1 we need to find the explicit formula for the function $F_{4t+3}(x)$ for $x \in [\frac{T_{4t+1}}{F_{4t+1}}, \frac{T_{4t}}{F_{4t}}]$.
\\

\textsc{Lemma 4.} \textit {For $x \in [\frac{T_{4t+1}}{F_{4t+1}}, \frac{T_{4t}}{F_{4t}}]$} 
\begin{equation}
\label{eq14}
F_{4t+3}(x)=\begin{cases}
F_{4t+1}x-T_{4t+1} &\text{if $x \in [\frac{T_{4t+1}}{F_{4t+1}}, \frac{T_{4t+3}+T_{4t+1}-1}{F_{4t+1}+F_{4t+3}}]$,}\\
T_{4t+3}-1-F_{4t+3}x &\text{if $x \in [\frac{T_{4t+3}+T_{4t+1}-1}{F_{4t+1}+F_{4t+3}}, \frac{T_{4t+3}-1}{F_{4t+3}}]$,}\\
F_{4t+3}x-T_{4t+3}+1 &\text{if $x \in [\frac{T_{4t+3}-1}{F_{4t+3}}, \frac{T_{4t+2}+T_{4t+3}-1}{F_{4t+3}+F_{4t+2}}]$,}\\
T_{4t+2}-F_{4t+2}x &\text{if $x \in [\frac{T_{4t+2}+T_{4t+3}-1}{F_{4t+3}+F_{4t+2}}, \frac{T_{4t+2}}{F_{4t+2}}]$,}\\
F_{4t+2}x-T_{4t+2} &\text{if $x \in [\frac{T_{4t+2}}{F_{4t+2}}, \frac{T_{4t+2}+1}{F_{4t+2}+3}]$}\\
1-3x &\text{if $x \in [\frac{T_{4t+2}+1}{F_{4t+2}+3}, \frac{T_{4t}-1}{F_{4t}-3}]$,}\\
T_{4t}-F_{4t}x &\text{if $x \in [\frac{T_{4t}-1}{F_{4t}-3}, \frac{T_{4t+3}-T_{4t}}{F_{4t+3}-F_{4t}}]$,}\\
T_{4t+3}-F_{4t+3}x &\text{if $x \in [\frac{T_{4t+3}-T_{4t}}{F_{4t+3}-F_{4t}}, \frac{T_{4t+3}}{F_{4t+3}}]$,}\\
F_{4t+3}x-T_{4t+3} &\text{if $x \in [\frac{T_{4t+3}}{F_{4t+3}}, \frac{T_{4t+3}+T_{4t}}{F_{4t+3}+F_{4t}}]$,}\\
T_{4t}-F_{4t}x &\text{if $x \in [\frac{T_{4t+3}+T_{4t}}{F_{4t+3}+F_{4t}}, \frac{T_{4t}}{F_{4t}}]$.}
\end{cases}
\end{equation}

Proof.

We prove Lemma 4 by induction.

For $t=1$ Lemma 4 follows from Lemma 1 (we are interested in the segment $[\frac{1}{5},\frac{1}{3}]$).

We assume that Lemma 4 is valid for $t=k$. 

To make our proof more clear we draw the graph of the function $F_{4k+3}(x)$ on the segment $[\frac{T_{4k+1}}{F_{4k+1}}, \frac{T_{4k}}{F_{4k}}]$ (see Fig.2). The domain below this graph is colored in black. The point $M$ is the intersection of the graphs of the functions $y=1-3x$ and $y=||F_{4k+2}x||$. The function $F_{4k+3}(x)$ attains its maximal value at this point. The scheme of the graph of the function $F_{4k+7}(x)$ on $[\frac{T_{4k+5}}{F_{4k+5}}, \frac{T_{4k+4}}{F_{4k+4}}]$ is marked with white.

The distance to the nearest integer from $\frac{F_{4k}}{\varphi+2}$, $\frac{F_{4k+1}}{\varphi+2}$, $\frac{F_{4k+2}}{\varphi+2}$, $\frac{F_{4k+3}}{\varphi+2}$ are also marked on the graph (according to (\ref{eq6}) - (\ref{eq9})).

\newpage

\begin{figure}[htbp]
\includegraphics[width=\textwidth]{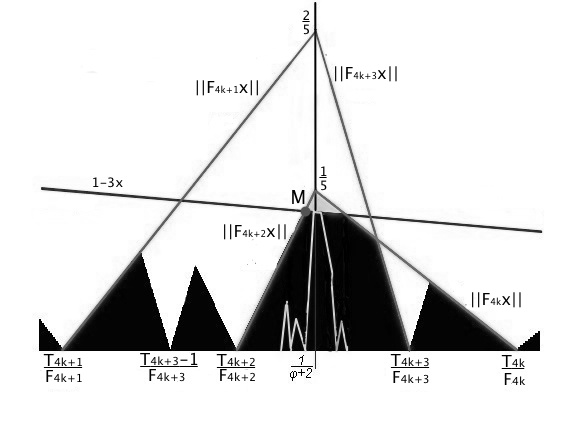}
\label{labelname}
\caption{$F_{4k+3}(x)$}
\end{figure}

\textsc{Remark 1.} \textit{To prove the inductive step we need to compare several numbers. These comparisons can be made in the following way:}

\textit{1) We have two expressions $A$ and $B$ depending on $T_{4k+7}$, $T_{4k+6}$, $T_{4k+5}$, $T_{4k+4}$, $T_{4k+3}$, $T_{4k+2}$, $F_{4k+7}$, $F_{4k+6}$, $F_{4k+5}$, $F_{4k+4}$, $F_{4k+3}$, $F_{4k+2}$, $F_{4k+1}$, $F_{4k}$. We want to prove that $A>B$. We consider $A-B$.}

\textit{2)  The values of $T_{i}$ are defined in (\ref{eq2})-(\ref{eq5}). We substitute these formulas into $A-B$. The expression obtained depends on $F_{4k+7}$, $F_{4k+6}$, $F_{4k+5}$, $F_{4k+4}$, $F_{4k+3}$, $F_{4k+2}$, $F_{4k+1}$, $F_{4k}$, $F_{4k-1}$.}

\textit{3) The following formulas can be obtained from the definition of Fibonacci sequence:}
\begin{center}

$F_{4k+7}=21F_{4k}+13F_{4k-1},$

$F_{4k+6}=13F_{4k}+8F_{4k-1},$

$F_{4k+5}=8F_{4k}+5F_{4k-1},$

$F_{4k+4}=5F_{4k}+3F_{4k-1},$

$F_{4k+3}=3F_{4k}+2F_{4k-1},$

$F_{4k+2}=2F_{4k}+F_{4k-1},$

$F_{4k+1}=F_{4k}+F_{4k-1}.$
\end{center}

\textit{We substitute them into the expression obtained in the second step. Now the difference $A-B$ depends on $F_{4k}$ and $F_{4k-1}$.}

\textit{4) We use the equality $(-1)^n=F_{n+1}F_{n-1}-F_{n}^2$ and the fact that $F_{4k}\ge3$ and $F_{4k-1}\ge2$ for $k\ge1$ to obtain inequality $A-B>0$.}
\\

Now we write the proof of one of such inequalities in the details. For example, $\frac{T_{4k+2}}{F_{4k+2}}<\frac{T_{4k+5}}{F_{4k+5}}$.

After the substitution described at the step 2 we obtain $$\frac{T_{4k+5}}{F_{4k+5}}-\frac{T_{4k+2}}{F_{4k+2}}=\frac{F_{4k+4}+2F_{4k+3}-2}{5F_{4k+5}}-\frac{3F_{4k}+F_{4k-1}-1}{5F_{4k+2}}.$$

After the substitution described at the step 3 we obtain $$\frac{T_{4k+5}}{F_{4k+5}}-\frac{T_{4k+2}}{F_{4k+2}}=\frac{11F_{4k}+7F_{4k-1}-2}{40F_{4k}+25F_{4k-1}}-\frac{3F_{4k}+F_{4k-1}-1}{10F_{4k}+5F_{4k-1}}.$$

We rewrite the expression obtained according to the step 4:

$$\frac{T_{4k+5}}{F_{4k+5}}-\frac{T_{4k+2}}{F_{4k+2}}=\frac{-10F^2_{4k}+10F_{4k-1}^2+10F_{4k}F_{4k-1}+20F_{4k}+15F_{4k-1}}{(40F_{4k}+25F_{4k-1})(10F_{4k}+5F_{4k-1})}=$$

$$=\frac{-10(F_{4k}F_{4k-2}-F_{4k-1}^2)+20F_{4k}+15F_{4k-1}}{(40F_{4k}+25F_{4k-1})(10F_{4k}+5F_{4k-1})}\ge\frac{-10(-1)^{4k-1}+60+30}{(40F_{4k}+25F_{4k-1})(10F_{4k}+5F_{4k-1})}>0.$$

The procedure described in Remark 1 will be used in our proof several times. Each time we use this procedure we refer to Remark 1.
\\

We move on to the proof of the inductive step. We assume that on the segment $[\frac{T_{4k+1}}{F_{4k+1}}, \frac{T_{4k}}{F_{4k}}]$ we know the explicit formula for $F_{4k+3}(x)$. We want to find the explicit formula for $F_{4k+7}(x)$ on the segment $[\frac{T_{4k+5}}{F_{4k+5}}, \frac{T_{4k+4}}{F_{4k+4}}]$.
\\

\textsc{Proposition 1.} \textit{The following inequalities are valid: }
$$\frac{T_{4k+2}}{F_{4k+2}}<\frac{T_{4k+5}}{F_{4k+5}}<\frac{T_{4k+7}-1}{F_{4k+7}}<\frac{T_{4k+6}}{F_{4k+6}}<\frac{1}{\varphi+2}<\frac{T_{4k+7}}{F_{4k+7}}<\frac{T_{4k+4}}{F_{4k+4}}<\frac{T_{4k+3}}{F_{4k+3}}.$$
\\

Proof.

To obtain the inequality $\frac{T_{4k+6}}{F_{4k+6}}<\frac{1}{\varphi+2}$ we should devide (\ref{eq8}) by $F_{4k+6}$. Similarly, the inequality $\frac{T_{4k+7}}{F_{4k+7}}>\frac{1}{\varphi+2}$ can be obtained from (\ref{eq9}).

The remaining inequalities of the proposition can be obtained by the means of procedure from Remark 1.

Proposition 1 is proved.
\\

\textsc{Proposition 2.} \textit{The zeros of the function $F_{4k+7}(x)$ on the segment $[\frac{T_{4k+5}}{F_{4k+5}}, \frac{T_{4k+4}}{F_{4k+4}}]$ are the points $\frac{T_{4k+5}}{F_{4k+5}}$, $\frac{T_{4k+7}-1}{F_{4k+7}}$, $\frac{T_{4k+6}}{F_{4k+6}}$, $\frac{T_{4k+7}}{F_{4k+7}}$, $\frac{T_{4k+4}}{F_{4k+4}}$.}
\\

Proof.

One can easily see that 
$$F_{4k+7}(x)=\min\{F_{4k+3}(x), ||F_{4k+4}x||, ||F_{4k+5}x||, ||F_{4k+6}x||, ||F_{4k+7}x||\}.$$
Hence we should find the zeros of functions $F_{4k+3}(x), ||F_{4k+4}x||, ||F_{4k+5}x||, ||F_{4k+6}x||, ||F_{4k+7}x||$ on the segment considered.

From the inductive assumption it follows that the function $F_{4k+3}(x)$ doesn't have zeros on $(\frac{T_{4k+2}}{F_{4k+2}}, \frac{T_{4k+3}}{F_{4k+3}})$. From Proposition 1 we know that $[\frac{T_{4k+5}}{F_{4k+5}}, \frac{T_{4k+4}}{F_{4k+4}}]\subset[\frac{T_{4k+2}}{F_{4k+2}}, \frac{T_{4k+3}}{F_{4k+3}}]$. So the function $F_{4k+3}(x)$ doesn't have zeros on the segment considered.
\\

The function $||F_{4k+4}x||$ is equal to zero only when $x=\frac{a}{F_{4k+4}}$, where $a$ is integer. The segment $[\frac{T_{4k+5}}{F_{4k+5}}, \frac{T_{4k+4}}{F_{4k+4}}]$ has only one point of such a kind which is $x=\frac{T_{4k+4}}{F_{4k+4}}$. We need to compare points $\frac{T_{4k+4}-1}{F_{4k+4}}$ and $\frac{T_{4k+4}+1}{F_{4k+4}}$ with the endpoints of the segment considered to prove that there are no other zeros on $[\frac{T_{4k+5}}{F_{4k+5}}, \frac{T_{4k+4}}{F_{4k+4}}]$.

So we need to prove that the inequalieties $\frac{T_{4k+4}-1}{F_{4k+4}}<\frac{T_{4k+5}}{F_{4k+5}}$ and $\frac{T_{4k+4}+1}{F_{4k+4}}>\frac{T_{4k+4}}{F_{4k+4}}$ are valid.

Obviously the inequality $\frac{T_{4k+4}+1}{F_{4k+4}}>\frac{T_{4k+4}}{F_{4k+4}}$ is valid. To prove the second one we need to use the procedure from Remark 1.
\\

The similar argument must be used to the analysis of the three remaining functions.

The function $||F_{4k+5}x||$ is equal to zero on the segment $[\frac{T_{4k+5}}{F_{4k+5}}, \frac{T_{4k+4}}{F_{4k+4}}]$ only when $x=\frac{T_{4k+5}}{F_{4k+5}}$. We need to prove that there are no other zeros of the function $||F_{4k+5}x||$ on the segment considered.

Obviously the inequlity $\frac{T_{4k+5}}{F_{4k+5}}>\frac{T_{4k+5}-1}{F_{4k+5}}$ is valid.

The inequality $\frac{T_{4k+5}+1}{F_{4k+5}}>\frac{T_{4k+4}}{F_{4k+4}}$ can be proved by the means of the procedure from Remark 1.
\\

The function $||F_{4k+6}x||$ is zero only when $x=\frac{T_{4k+6}}{F_{4k+6}}$. From Proposition 1 it follows that $\frac{T_{4k+6}}{F_{4k+6}}\in[\frac{T_{4k+5}}{F_{4k+5}}, \frac{T_{4k+4}}{F_{4k+4}}]$. We need to prove that $||F_{4k+6}x||$ doesn't have other zeros on the segment considered.

The inequalites $\frac{T_{4k+5}}{F_{4k+5}}>\frac{T_{4k+6}-1}{F_{4k+6}}$ and $\frac{T_{4k+6}+1}{F_{4k+6}}>\frac{T_{4k+4}}{F_{4k+4}}$ are proved by the means of the procedure from Remark 1.
\\

The function $||F_{4k+7}x||$ is zero when $x=\frac{T_{4k+7}-1}{F_{4k+7}}$ and $x=\frac{T_{4k+7}}{F_{4k+7}}$. From Proposition 1 we know that these points belong to the segment $[\frac{T_{4k+5}}{F_{4k+5}}, \frac{T_{4k+4}}{F_{4k+4}}]$. We show that the function $||F_{4k+7}x||$ doesn't have other zeros.

By the means of the procedure from Remark 1 we obtain the inequalities $\frac{T_{4k+5}}{F_{4k+5}}>\frac{T_{4k+7}-2}{F_{4k+7}}$ and $\frac{T_{4k+7}+1}{F_{4k+7}}>\frac{T_{4k+4}}{F_{4k+4}}$.

Proposition 2 is proved.
\\

\textsc{Proposition 3.} \textit{The equalities $\frac{T_{4k+6}}{F_{4k+6}}=\frac{T_{4k+2}+1}{F_{4k+2}+3}$ and $\frac{T_{4k+4}}{F_{4k+4}}=\frac{T_{4k}-1}{F_{4k}-3}$ are valid. And}
\begin{equation}
\label{eq15}
F_{4k+3}(x)=\begin{cases}
F_{4k+2}x-T_{4k+2} &\text{if $x \in [\frac{T_{4k+5}}{F_{4k+5}}, \frac{T_{4k+6}}{F_{4k+6}}]$,}\\
1-3x &\text{if $x \in [\frac{T_{4k+6}}{F_{4k+6}}, \frac{T_{4k+4}}{F_{4k+4}}]$.}
\end{cases}
\end{equation}

Proof.

We prove the equality $\frac{T_{4k+6}}{F_{4k+6}}=\frac{T_{4k+2}+1}{F_{4k+2}+3}$:

$$
\frac{T_{4k+6}}{F_{4k+6}}-\frac{T_{4k+2}+1}{F_{4k+2}+3}=\frac{3F_{4k+4}+F_{4k+3}-1}{5F_{4k+6}}-\frac{3F_{4k}+F_{4k-1}+4}{5F_{4k+2}+15}=
$$

$$
=\frac{18F_{4k}+11F_{4k-1}-1}{65F_{4k}+40F_{4k-1}}-\frac{3F_{4k}+F_{4k-1}+4}{10F_{4k}+5F_{4k-1}+15}=
$$

$$
=\frac{-15(F_{4k}F_{4k-2}-F_{4k-1}^2)-15}{(65F_{4k}+40F_{4k-1})(10F_{4k}+5F_{4k-1}+15)}=0
$$

This equality has the following meaning: the finction $F_{4t+3}(x)$ attains its maximal value just at the zero of the function $F_{4t+7}$.
\\

We prove the equality $\frac{T_{4k+4}}{F_{4k+4}}=\frac{T_{4k}-1}{F_{4k}-3}$:

$$
\frac{T_{4k}-1}{F_{4k}-3}-\frac{T_{4k+4}}{F_{4k+4}}=\frac{2F_{4k}-F_{4k-1}-4}{5(F_{4k}-3)}-\frac{2F_{4k+4}-F_{4k+3}+1}{5F_{4k+4}}=
$$

$$
=\frac{2F_{4k}-F_{4k-1}-4}{5F_{4k}-15}-\frac{7F_{4k}+4F_{4k-1}+1}{25F_{4k}+15F_{4k-1}}=
$$

$$
=\frac{15(F_{4k}F_{4k-2}-F_{4k-1}^2)+15}{(5F_{4k}-15)(25F_{4k}+15F_{4k-1})}=0
$$

These two equalities and the inequality $\frac{T_{4k+2}}{F_{4k+2}}<\frac{T_{4k+5}}{F_{4k+5}}$ from Proposition 1 lead to formula (\ref{eq15}).

Proposition 3 is proved.
\\

\textsc{Proposition 4.} \textit{The following equalities are valid:}

$$
||F_{4k+4}x||=
T_{4k+4}-F_{4k+4}x \text{\quad if $x \in [\frac{T_{4k+5}}{F_{4k+5}}, \frac{T_{4k+4}}{F_{4k+4}}]$,}
$$

$$
||F_{4k+5}x||=\begin{cases}
F_{4k+5}x-T_{4k+5} &\text{if $x\in[\frac{T_{4k+5}}{F_{4k+5}}, \frac{2T_{4k+5}+1}{2F_{4k+5}}]$,}\\
T_{4k+5}+1-F_{4k+5}x &\text{if $x \in [\frac{2T_{4k+5}+1}{2F_{4k+5}}, \frac{T_{4k+4}}{F_{4k+4}}]$,}
\end{cases}
$$

$$
||F_{4k+6}x||=\begin{cases}
T_{4k+6}-F_{4k+6}x &\text{if $x\in[\frac{T_{4k+5}}{F_{4k+5}}, \frac{T_{4k+6}}{F_{4k+6}}]$,}\\
F_{4k+6}x-T_{4k+6} &\text{if $x \in [\frac{T_{4k+6}}{F_{4k+6}}, \frac{2T_{4k+6}+1}{2F_{4k+6}}]$,}\\
T_{4k+6}+1-F_{4k+6}x &\text{if $x \in[\frac{2T_{4k+6}+1}{2F_{4k+6}}, \frac{T_{4k+4}}{F_{4k+4}}]$,}
\end{cases}
$$

$$
||F_{4k+7}x||=\begin{cases}
T_{4k+7}-1-F_{4k+7}x &\text{if $x\in[\frac{T_{4k+5}}{F_{4k+5}}, \frac{T_{4k+7}-1}{F_{4k+7}}]$,}\\
F_{4k+7}x-T_{4k+7}+1 &\text{if $x \in [\frac{T_{4k+7}-1}{F_{4k+7}}, \frac{2T_{4k+7}-1}{2F_{4k+7}}]$,}\\
T_{4k+7}-F_{4k+7}x &\text{if $x \in[\frac{2T_{4k+7}-1}{2F_{4k+7}}, \frac{T_{4k+7}}{F_{4k+7}}]$,}\\
F_{4k+7}x-T_{4k+7} &\text{if $x \in [\frac{T_{4k+7}}{F_{4k+7}}, \frac{T_{4k+4}}{F_{4k+4}}]$.}
\end{cases}
$$

Proof.

For the explicit formula of the function $||F_nx||$ we shoud know the zeros of this function as well as the points where this function attains its maximal value.

From Proposition 2 we know the zeros of the functions $||F_{4k+4}x||$, $||F_{4k+5}x||$, $||F_{4k+6}x||$, $||F_{4k+7}x||$ on the segment considered. The maximal value of the function $||F_nx||$ is attained at the middle between two neighbouring zeros of this function. So it's enough to prove the following inequlities:

1) $\frac{T_{4k+5}}{F_{4k+5}}>\frac{2T_{4k+4}-1}{2F_{4k+4}}$ (for the middle of the segment $[\frac{T_{4k+4}-1}{F_{4k+4}}, \frac{T_{4k+4}}{F_{4k+4}}]$) 

2) $\frac{2T_{4k+5}+1}{2F_{4k+5}}<\frac{T_{4k+4}}{F_{4k+4}}$ (for the middle of the segment $[\frac{T_{4k+5}}{F_{4k+5}}, \frac{T_{4k+5}+1}{F_{4k+5}}]$),

3) $\frac{2T_{4k+6}-1}{2F_{4k+6}}<\frac{T_{4k+5}}{F_{4k+5}}$ (for the middle of the segment $[\frac{T_{4k+6}-1}{F_{4k+6}}, \frac{T_{4k+6}}{F_{4k+6}}]$), 

4) $\frac{2T_{4k+6}+1}{2F_{4k+6}}<\frac{T_{4k+4}}{F_{4k+4}}$ (for the middle of the segment $[\frac{T_{4k+6}}{F_{4k+6}}, \frac{T_{4k+6}+1}{F_{4k+6}}]$), 

5) $\frac{2T_{4k+7}-3}{2F_{4k+7}}<\frac{T_{4k+5}}{F_{4k+5}}$ (for the middle of the segment $[\frac{T_{4k+7}-2}{F_{4k+7}}, \frac{T_{4k+7}-1}{F_{4k+7}}]$), 

6) $\frac{2T_{4k+7}+1}{2F_{4k+7}}>\frac{T_{4k+4}}{F_{4k+4}}$ (for the middle of the segment $[\frac{T_{4k+7}}{F_{4k+7}}, \frac{T_{4k+7}+1}{F_{4k+7}}]$).

All these inequalities are proved by the means of the procedure from Remark 1.

These inequlities mean that the points $\frac{2T_{4k+5}+1}{2F_{4k+5}}$, $\frac{2T_{4k+6}+1}{2F_{4k+6}}$ belong to the segment $[\frac{T_{4k+5}}{F_{4k+5}}, \frac{T_{4k+6}}{F_{4k+6}}]$, the others considered do not.

Proposition 4 is proved.
\\

Now we know the explicit formula for the functions $F_{4k+3}(x)$, $||F_{4k+4}x||$, $||F_{4k+5}x||$, $||F_{4k+6}x||$, $||F_{4k+7}x||$ on the segment $[\frac{F_{4k+5}}{T_{4k+5}}, \frac{T_{4k+4}}{F_{4k+4}}]$ and the relative position of zeros of these functions on this segment (from Proposition 1). So the way to find the explicit formula for the function $F_{4k+7}(x)$ on $[\frac{F_{4k+5}}{T_{4k+5}}, \frac{T_{4k+4}}{F_{4k+4}}]$ become obvious. By the means of the procedure from Remark 1 for comparison of two numbers the formula (\ref{eq14}) is obtained.

Lemma 4 is proved.
\\

\textsc{Lemma 5.} \textit {The graph of the function $F_{4t+3}(x)$ on the segment $[0,\frac{1}{2}]$ has only one vertex which lies above the line $y=\frac{\varphi-1}{\varphi+2}$. This vertex has coordinates $x=\frac{T_{4t+2}+1}{F_{4t+2}+3}$, $y=1-3\frac{T_{4t+2}+1}{F_{4t+2}+3}$. Moreover we have $x \in [\frac{T_{4t+2}}{F_{4t+2}}, \frac{T_{4t+3}}{F_{4t+3}}]$.}
\\

Proof.

We prove Lemma 5 by induction.

For $t=1$ the statement of lemma is obtained from Lemma 2.

We assume that for $t=k$ Lemma 5 is valid. 

The point $(\frac{T_{4k+2}+1}{F_{4k+2}}, 1-3\frac{T_{4k+2}}{F_{4k+2}})$ is the intersection of the lines $y=F_{4k+2}x-T_{4k+2}$ and $y=1-3x$. These lines intersect the line $y=\frac{\varphi-1}{\varphi+2}$ in points $x=\frac{\varphi-1+T_{4k+2}(\varphi+2)}{(\varphi+2)F_{4k+2}}$ and $x=\frac{1}{\varphi+2}$ correspondently. So from the assumption of induction it follows that $F_{4k+3}(x)\ge\frac{\varphi-1}{\varphi+2}$ only if $x \in [\frac{\varphi-1+T_{4k+2}(\varphi+2)}{(\varphi+2)F_{4k+2}}, \frac{1}{\varphi+2}]$.

We prove that the statement of lemma is valid for $t=k+1$. 

We note that the points $x=\frac{\varphi-1+T_{4k+2}(\varphi+2)}{(\varphi+2)F_{4k+2}}$ and $x=\frac{1}{\varphi+2}$ belong to the segment $[\frac{T_{4k+5}}{F_{4k+5}}, \frac{T_{4k+4}}{F_{4k+4}}]$. Indeed, for the point $\frac{1}{\varphi+2}$ it follows from Proposition 1. For the other one we prove the inequality $\frac{T_{4k+5}}{F_{4k+5}}<\frac{\varphi-1+T_{4k+2}(\varphi+2)}{(\varphi+2)F_{4k+2}}$ by the means of the procedure from Remark 1 and the bound $1,618<\varphi<1,619$.

So $F_{4k+3}(x)\le\frac{\varphi-1}{\varphi+2}$ outside $[\frac{T_{4k+5}}{F_{4k+5}}, \frac{T_{4k+4}}{F_{4k+4}}]$. Since $F_{4k+7}(x)\le F_{4k+3}(x)$, then $F_{4k+7}(x)\le \frac{\varphi-1}{\varphi+2}$ outside this segment.

From Lemma 4 we know the explicit formula for $F_{4k+7}(x)$ on the segment $[\frac{T_{4k+5}}{F_{4k+5}}, \frac{T_{4k+4}}{F_{4k+4}}]$. This formula leads to the fact that the only vertex of the graph is above the line $y=\frac{\varphi-1}{\varphi+2}$. It has the coordinates $(\frac{T_{4k+6}+1}{F_{4k+6}+3}, 1-3\frac{T_{4k+2}+1}{F_{4k+2}+3}$).

Lemma 5 is proved.
\\

So we see that the maximum of the function $F_{4t+3}(x)$ is equal to $1-3\frac{T_{4t+2}+1}{F_{4t+2}+3}$. This function attains its maximal value at $x=\frac{T_{4t+2}+1}{F_{4t+2}+3}$.
\\

\textsc{Proposition 5.} \textit{The equalities $\frac{T_{4t+2}+1}{F_{4t+2}+3}=\frac{F_{2t+2}}{F_{2t+2}+F_{2t+4}}$ and $1-3\frac{T_{4t+2}+1}{F_{4t+2}+3}=\frac{F_{2t+1}}{F_{2t+2}+F_{2t+4}}$ are valid.}
\\

Proof.

To prove the first equality it's enough to show that the following equlities are valid: 

$$F_{4t+2}+3=F_{2t-1}(F_{2t+2}+F_{2t+4}),\qquad T_{4t+2}+1=F_{2t-1}F_{2t+2}.$$

We prove the equality $F_{4t+2}+3-F_{2t-1}(F_{2t+2}+F_{2t+4})=0$. As for Fibonacci numbers we have the formula 
$F_{2t}=F_{t+1}^2-F_{t-1}^2$
we deduce
$$
F_{4t+2}+3-F_{2t-1}(F_{2t+2}+F_{2t+4})=F_{2t+2}^2-F_{2t}^2+3-F_{2t-1}F_{2t+2}-F_{2t-1}F_{2t+4}.
$$
Then we substitute the formulas 
$$
\begin{cases}
F_{2t-1}=F_{2t+1}-F_{2t},\\
F_{2t+2}=F_{2t+1}+F_{2t},\\
F_{2t+4}=3F_{2t+1}+2F_{2t}.
\end{cases}
$$
into the expression obtained.
This substitution leads to the following equalities:
$$
F_{2t+2}^2-F_{2t}^2+3-F_{2t-1}F_{2t+2}-F_{2t-1}F_{2t+4}=3-3F_{2t+1}^2+3F_{2t}^2+3F_{2t+1}F_{2t}=
$$
$$
=3-3F_{2t+1}F_{2t-1}+3F_{2t}^2=3-(-3)^{2t}=0.
$$

The equality $T_{4t+2}+1=F_{2t-1}F_{2t+2}$ is proved in the same way. 
\\

Since the equality $\frac{T_{4t+2}+1}{F_{4t+2}+3}=\frac{F_{2t+2}}{F_{2t+2}+F_{2t+4}}$ is valid, then the proof of the equality $1-3\frac{T_{4t+2}+1}{F_{4t+2}+3}=\frac{F_{2t+1}}{F_{2t+2}+F_{2t+4}}$ is obvious:

$$1-3\frac{T_{4t+2}+1}{F_{4t+2}+3}=1-3\frac{F_{2t+2}}{F_{2t+2}+F_{2t+4}}=\frac{F_{2t+2}+F_{2t+4}-3F_{2t+2}}{F_{2t+2}+F_{2t+4}}=\frac{F_{2t+1}}{F_{2t+2}+F_{2t+4}}$$

Proposition 5 is proved.
\\

Theorem 1 is proved.

\section{Proof of Theorem 2}

From Lemma 1, Lemma 4 and Proposition 1 we see that
$$\min_{n=1,\ldots,N}||{F_{n}}\alpha_1||=1-3\alpha_1=\frac{\varphi-1}{\varphi+2}.$$

Hence Theorem 2 is proved.

\section{Proof of Theorem 3}

From (\ref{eq6}) - (\ref{eq9}) we see that $||\frac{F_{4t}}{\varphi+2}||=\frac{1}{5}+O(\varphi^{-4t})$, $||\frac{F_{4t+1}}{\varphi+2}||=-\frac{2}{5}+O(\varphi^{-4t})$, $||\frac{F_{4t+2}}{\varphi+2}||=-\frac{1}{5}+O(\varphi^{-4t})$, $||\frac{F_{4t}}{\varphi+2}||=\frac{2}{5}+O(\varphi^{-4t})$. Therefore Theorem 3 is proved.

\section{Proof of Theorem 4}

Put $G_1(x,y)=x$, $G_2(x,y)=y$. Then $G_n (x,y)=G_{n-1}+G_{n-2}$ for $n\ge3$. Let $G^N(x,y)=\min_{n=1 \ldots N} ||G_n(x,y)||$. Put $t_N=\max_{x,y \in \mathbb R} G^N(x,y)$.
 
For each $N$ we devide the square $[0,1]\times[0,1]$ into domains with the same minimal functions.

The maximum is attained at the boundaries of the neighbouring domains.

The graphs below represent such division into domains for $N=2,...5$. For the first three graphs the minimal distance to the nearest integer is indicated in each domain. The points where $G^N(x,y)$ attains its maximum, the boundaries of the domains, the lines where $G^N(x,y)=0$ are also marked.

\begin{multicols}{2}
\includegraphics[width=50mm]{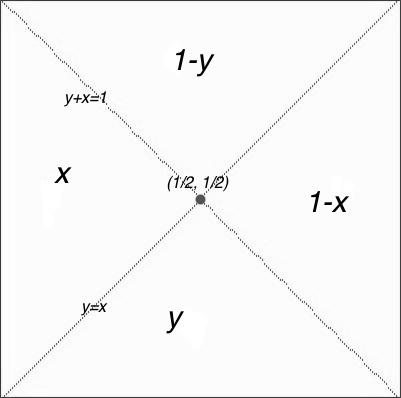}
\\
$\text{\qquad \qquad \qquad Figure 3}$

For $N=2$ (see Fig.3)
$t_N=\frac{1}{2}$, 

which attains at $x=\frac{1}{2}$, $y=\frac{1}{2}$.
\\

\includegraphics[width=50mm]{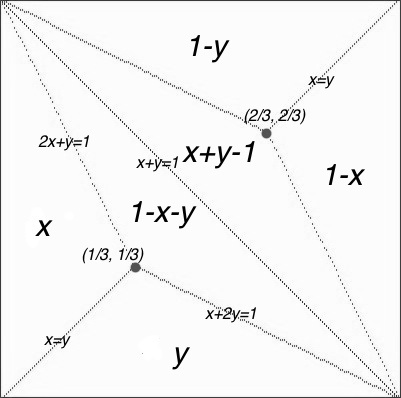}
\\
$\text{\qquad \qquad \qquad Figure 4}$
\\

For $N=3$ (see Fig.4) $t_N=\frac{1}{3}$, 

which attains at:

$x_1=\frac{1}{3}$, $y_1=\frac{1}{3}$, $||G_3(x_1,y_1)||=\frac{2}{3}$,

$x_2=\frac{2}{3}$, $y_2=\frac{2}{3}$, $||G_3(x_2,y_2)||=\frac{1}{3}$.
\end{multicols}

\includegraphics[width=50mm]{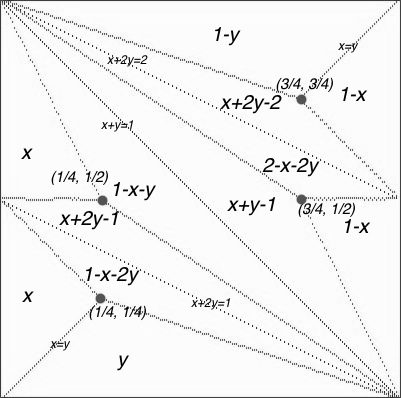}
\\$\text{\qquad \qquad \qquad Figure 5}$
\\

For $N=4$ (see Fig. 5) $t_N=\frac{1}{4}$, 

which attains at:

$x_1=\frac{1}{4}$, $y_1=\frac{1}{4}$, $||G_3(x_1, y_1)||=\frac{2}{4}$, $||G_4(x_1, y_1)||=\frac{3}{4}$,

$x_2=\frac{1}{4}$, $y_2=\frac{2}{4}$, $||G_3(x_2, y_2)||=\frac{3}{4}$, $||G_4(x_2, y_2)||=\frac{1}{4}$,

$x_3=\frac{3}{4}$, $y_3=\frac{3}{4}$, $||G_3(x_3, y_3)||=\frac{2}{4}$, $||G_4(x_3, y_3)||=\frac{1}{4}$,

$x_4=\frac{3}{4}$, $y_4=\frac{2}{4}$, $||G_3(x_4, y_4)||=\frac{1}{4}$, $||G_4(x_4, y_4)||=\frac{3}{4}$.
\\

On the next two graphs we only mark the lines where $G^N(x,y)$ is equal to zero and the points where this function attains its maximal value.

\includegraphics[width=\textwidth]{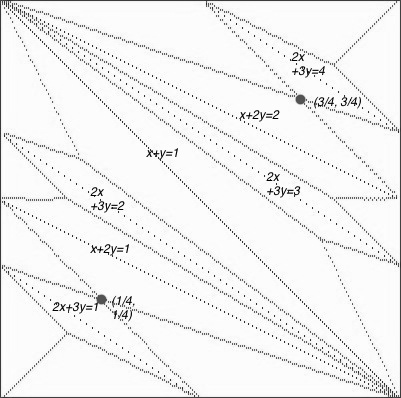}
\\$\text{\qquad \qquad \qquad  \qquad \qquad
\qquad \qquad \qquad \qquad \qquad Figure 6}$
\\

For $N=5$ (see Fig. 6) $t_5=\frac{1}{4}$, which attains at:

$x_1=\frac{1}{4}$, $y_1=\frac{1}{4}$, $||G_3(x_1, y_1)||=\frac{2}{4}$, $||G_4(x_1, y_1)||=\frac{3}{4}$, $||G_5(x_1, y_1)||=\frac{1}{4},$

$x_2=\frac{3}{4}$, $y_2=\frac{3}{4}$, $||G_3(x_2, y_2)||=\frac{2}{4}$, $||G_4(x_2, y_2)||=\frac{1}{4}$, $||G_5(x_2, y_2)||=\frac{3}{4}.$
\\

\newpage

\includegraphics[width=\textwidth]{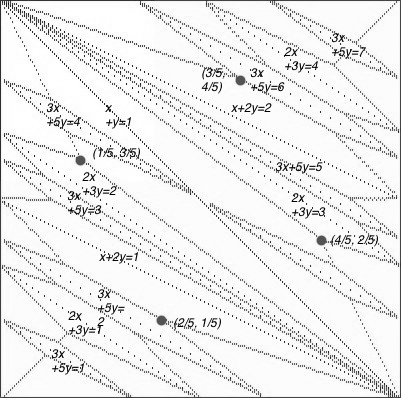}
\\$\text{\qquad \qquad \qquad  \qquad \qquad
\qquad \qquad \qquad \qquad \qquad Figure 7}$
\\

For $N=6$ (see Fig. 7) $t_6=\frac{1}{5}$, which attains at:

$x_1=\frac{2}{5}$,$y_1=\frac{1}{5}$,$||G_3(x_1,y_1)||=\frac{3}{5}$,$||G_4(x_1, y_1)||=\frac{4}{5}$,$||G_5(x_1, y_1)||=\frac{2}{5}$,$||G_6(x_1, y_1)||=\frac{1}{5},$
\\

$x_2=\frac{1}{5}$,$y_2=\frac{3}{5}$,$||G_3(x_2,y_2)||=\frac{4}{5}$,$||G_4(x_2,y_2)||=\frac{2}{5}$,$||G_5(x_2, y_2)||=\frac{1}{5}$,$||G_6(x_2,y_2)||=\frac{3}{5},$
\\

$x_3=\frac{3}{5}$,$y_3=\frac{4}{5}$,$||G_3(x_3,y_3)||=\frac{2}{5}$,$||G_4(x_3,y_3)||=\frac{1}{5}$,$||G_5(x_3,y_3)||=\frac{3}{5}$,$||G_6(x_3,y_3)||=\frac{4}{5},$
\\

$x_4=\frac{4}{5}$,$y_4=\frac{2}{5}$,$||G_3(x_4,y_4)||=\frac{1}{5}$,$||G_4(x_4,y_4)||=\frac{3}{5}$,$||G_5(x_4,y_4)||=\frac{4}{5}$,$||G_6(x_4,y_4)||=\frac{2}{5}$.
\\

We note, that each of these sequences is periodic. For example, for $x=\frac{2}{5}$ and $y=\frac{1}{5}$ the sequence $G_N(x,y)$ is as follows: $\frac{2}{5}$, $\frac{1}{5}$, $\frac{3}{5}$, $\frac{4}{5}$, $\frac{2}{5}$, $\frac{1}{5}$, $\frac{3}{5}$, $\frac{4}{5}$, $\frac{2}{5}$, $\frac{1}{5}\ldots$. It means that $||G^N(\frac{2}{5}, \frac{1}{5})||=\frac{1}{5}$ for $N\ge2$.
\\

We note that enequality $G^{N+1}(x,y) \le G^{N}(x,y)$ is always valid. Hence, $t_N=\frac{1}{5}$ for $N\ge6$ .

Theorem 4 is proved.

\textit{Victoria Zhuravleva}\footnote{Research is supported by RFBR Grant 12-01-00681}

\textit{Moscow Lomonosov State University}

\textit{e-mail: v.v.zhuravleva@gmail.com}


\begin{thebibliography}{xx}

\bibitem{ref1} \textsc{R.K. Akhunzhanov}, 
\textit{On the distribution modulo 1 of exponential sequences}.
Mathematical notes {\bf76:2} (2004), 153--160.


\bibitem{ref3} {\sc A. Dubickas}, 
\textit{Arithmetical properties of powers of algebraic numbers}. 
Bull. London Math. Soc. {\bf38} (2006), 70--80.

\bibitem{ref4} {\sc L. Kuipers, H. Niederreiter}, 
\textit{Uniform distribution of sequences}. 
John Wiley \& Sons (1974).

\bibitem{ref5} {\sc W.M. Schmidt}, 
\textit{Diophantine approximations}. 
Lect. Not. Math. {\bf785} (1980).

\bibitem{ref6} {\sc W.M. Schmidt}, 
\textit{On badly approximable numbers and certain games}. 
Trans. Amer. Math. Soc. {\bf623} (1966), 178--199.
\end{thebibliography}
\end{document}